\documentstyle[12pt]{article}

 2
\def\BI{{\rm 1\!l}}
\def\pois#1#2{\{#1,#2\}}

\def\Eq#1{{\begin{equation} #1 \end{equation}}}

\def\gone{\matrix{{}\cr g \cr {}^1}}
\def\gtwo{\matrix{{}\cr g \cr {}^2}}

\def\gbone{\matrix{{}\cr \bar g \cr {}^1}}
\def\gbtwo{\matrix{{}\cr \bar g \cr {}^2}}
\def\jone{\matrix{{}\cr j \cr {}^1}}
\def\jtwo{\matrix{{}\cr j \cr {}^2}}

\def\jdtwo{\matrix{{}\cr j^\dagger \cr {}^2}}
\def\xone{\matrix{{}\cr x \cr {}^1}}
\def\xtwo{\matrix{{}\cr x \cr {}^2}}
\def\stwo{\matrix{{}\cr s \cr {}^2}}

\def\ptone{\matrix{{}\cr \tilde p \cr {}^1}}
\def\pttwo{\matrix{{}\cr \tilde p \cr {}^2}}
\def\gamone{\matrix{{}\cr  \gamma \cr {}^1}}
\def\gamtwo{\matrix{{}\cr  \gamma \cr {}^2}}
\def\gambone{\matrix{{}\cr  \bar\gamma \cr {}^1}}
\def\gambtwo{\matrix{{}\cr  \bar\gamma \cr {}^2}}
\def\wtone{\matrix{{}\cr \tilde w \cr {}^1}}
\def\wttwo{\matrix{{}\cr \tilde w \cr {}^2}}
\def\fone{\matrix{{}\cr f \cr {}^1}}
\def\ftwo{\matrix{{}\cr f \cr {}^2}}
\def\fdtwo{\matrix{{}\cr f^\dagger \cr {}^2}}
\def\Jtwo{\matrix{{}\cr  J \cr {}^2}}
\def\Jdone{\matrix{{}\cr  J^\dagger \cr {}^1}}
\def\sJtwo{\sin (\lambda \Jtwo)}
\def\sJdone{\sin (\lambda \Jdone)}
\def\cJdone{\cos (\lambda \Jdone)}
\def\c2Jtwo{\cos^2 (\lambda \Jtwo)}
\def\cJtwo{\cos (\lambda \Jtwo)}
\def\Tone{\matrix{{}\cr T \cr {}^1}}
\def\Ttwo{\matrix{{}\cr T \cr {}^2}}
\def\Tbone{\matrix{{}\cr \bar T \cr {}^1}}
\def\Tbtwo{\matrix{{}\cr \bar T \cr {}^2}}
\def\Gone{\matrix{{}\cr \Gamma \cr {}^1}}
\def\Gtwo{\matrix{{}\cr \Gamma \cr {}^2}}
\def\Gbtwo{\matrix{{}\cr \bar \Gamma \cr {}^2}}
\def\Gbone{\matrix{{}\cr \bar \Gamma \cr {}^1}}
\def\Rontw{\matrix{{}\cr R \cr {}^{12}}}
\def\Rtwon{\matrix{{}\cr R \cr {}^{21}}}

\def\Pone{\matrix{{}\cr P \cr {}^1}}
\def\Ptwo{\matrix{{}\cr P \cr {}^2}}
\begin{document}
\begin{flushright}
UAHEP 952\\
May 1995\\
\end{flushright}
%\vspace*{5mm}

\centerline{ \LARGE  Lie-Poisson Deformation of the Poincar\'e Algebra}
\vskip 2cm

\centerline{ {\sc  A. Stern${}^{*\dagger}$ and I. Yakushin${}^*$ }  }

\vskip 1cm

\centerline{ ${}^{*}$ Dept. of Physics and Astronomy, Univ. of Alabama,
Tuscaloosa, Al 35487, U.S.A.}

\vskip 1cm

\centerline{ ${}^{\dagger}$ The Erwin Schr\"odinger International
Institute for Mathematical Physics,}
\centerline{Pasteurgasse 6/7, A-1090 Wien, Austria}

\vskip 2cm

\vspace*{5mm}

\normalsize
\centerline{\bf ABSTRACT}
We find a one parameter family of quadratic
Poisson structures on ${\bf R}^4\times SL(2,C)$  which
satisfies the property {\it a)} that it is preserved under the Lie-Poisson
action of the Lorentz group, as well as {\it b)}
that it reduces to the standard
Poincar\'e algebra for a particular limiting value of the
parameter.   (The Lie-Poisson
transformations reduce to canonical ones in that limit, which we
therefore refer to as the `canonical limit'.)
Like with the Poincar\'e algebra, our
 deformed Poincar\'e algebra has two Casimir functions which we associate
with `mass' and `spin'.  We parametrize the symplectic leaves
of ${\bf R}^4\times SL(2,C)$ with space-time coordinates, momenta
and spin, thereby  obtaining realizations of the deformed
algebra for the cases of a spinless and a spinning particle.
The formalism can be applied for finding a one
parameter family of canonically inequivalent descriptions of
the photon.

\vskip 2cm
\vspace*{5mm}

\newpage
\scrollmode
\section{Introduction}

A number of authors have examined
deformations of the Poincar\'e algebra in quantum theory,
in addition to investigating the effects  of deforming
the usual space-time symmetry.[1-7]  Of course there exists no
unique procedure for carrying out such deformations.
The proposals made so far are generally
within the framework of Hopf algebras, and they often rely upon
making a  contraction of the quantum de Sitter algebra.

Another technique is to deform the Poincar\'e algebra already
at the classical level.\cite{Zak}  There the deformed
 algebra is to be realized
in terms of Poisson brackets rather than commutation
relations, and the construction should be made within the framework
 of Poisson-Lie groups in order
to make connection with Hopf algebras in quantum theory.
This is the approach we shall follow here.

The classical analysis is considerably simpler
than its quantum counterpart for a number of reasons.
One reason is of courses that the elements of the algebra, i.e.
the classical observables, are
commuting variables.  In addition to this, to check the
consistency of the algebra we essentially only need to verify the Jacobi
identity (although this may not always be so easy).
 Furthermore, in the classical theory
symmetries are associated with ordinary Lie groups, and not
quantum groups.

With regard to the classical symmetries, which we denote by ${\cal
S}$, we shall here only be    concerned with
the Lorentz group [or actually, $SL(2,C)$], which we here regard as
a Poisson-Lie group.  Thus unlike in canonical theories, the symmetry group
carries a Poisson structure, which by definition is consistent
 with left or right group multiplication.   Such
structures are well known.  Our choice of Poisson brackets
$\{,\}_{\cal S}$ corresponds to the classical
analogue of the defining relations for  $SL_q(2,C)$.\cite{Tak}

With regard to the space of classical observables, which we denote by ${\cal
O}$, we examine    a one parameter family of algebras defined on
${\bf R}^4\times SL(2,C)$.   These algebras
are preserved under the Lie-Poisson action of ${\cal S}$.  Furthermore,
they  reduce to
the Poincar\'e algebra for a particular limiting value of the
parameter.  We refer to the limiting value as the `canonical limit'.
This is because ${\cal S}$ has a trivial Poisson bracket in the limit,
i.e. $\{,\}_{\cal S}\rightarrow 0$, and hence the Lie-Poisson action
of ${\cal S}$ on ${\cal O}$ simplifies to the canonical action in  the  limit.

For all values of the deformation parameter,
${\cal S}$ will act on ${\cal O}$ in the
standard way, i.e. momenta transforms as a Lorentz vector, while
angular momenta transforms with the adjoint representation of
$SL(2,C)$.  This defines a map,
 $\sigma:{\cal S}\times {\cal O}\rightarrow {\cal O}$.
Because ${\cal S}$ induces the Lie-Poisson  action on ${\cal O}$,
 $\sigma$ must be a Poisson map, which means that if $f_1$ and $f_2$ are
functions on  ${\cal O}$, then
 \Eq{  \{f_1,f_2\}_{\cal O} \circ\sigma =
 \{ f_1\circ\sigma, f_2\circ\sigma\}_{{\cal O}\times {\cal S}} \;,}
where the product Poisson structure is assumed on
${\cal O}\times {\cal S}$, which means that the symmetries have zero
Poisson brackets with the classical observables.

Our  deformed Poincar\'e   algebra on ${\cal
O}$ can be  completely specified by four quadratic
Poisson bracket relations $\{,\}_{\cal O}$ which we give below.
The brackets are evaluated between
variables which we shall associate with
`momenta' and `angular momenta' spanning ${\bf R}^4\times SL(2,C)$.
 The momenta will be expressed in terms
of a $2\times 2$ hermitean matrix $\tilde p$, while the angular momenta
are contained in a $2\times 2$ complex unimodular matrix $\gamma$.
The four Poisson bracket relations
 can be written in terms of a classical $r-$matrix
(and its hermitean conjugate $r^\dagger$) which is assumed to
 satisfy the (modified) classical Yang-Baxter equations.
Using tensor product notation, the four relations are:
\begin{eqnarray}
\pois{\ptone}{\pttwo} &= & r \ptone\pttwo +\ptone\pttwo
r^\dagger  -\pttwo r^\dagger \ptone -\ptone r \pttwo \;,\label{pbptpt}\\
\pois{\gamone}{\gamtwo}& = & r^\dagger \gamone\gamtwo +\gamone\gamtwo
r  -\gamtwo r \gamone -\gamone r^\dagger \gamtwo \;,\label{pbgam}\\
 \{\gamone,\gambtwo\}&=&
  r \gamone\gambtwo +\gamone\gambtwo r
  -\gambtwo r \gamone -\gamone r \gambtwo \;,\label{pbgmgmb}  \\
 \{\ptone,\gamtwo\}&=&
  r^\dagger \ptone\gamtwo +\ptone\gamtwo r
  -\gamtwo r^\dagger \ptone -\ptone r^\dagger \gamtwo \;,\label{pbpgm}
\end{eqnarray}
where $\bar\gamma =  {\gamma^\dagger}^{-1}$.
The $1$ and $2$ labels refer to two separate vector spaces, with
$\ptone=\tilde p\otimes \BI,\;\pttwo=\BI\otimes \tilde p ,$
 $\gamone=\gamma\otimes \BI$ and $  \gamtwo=\BI\otimes \gamma ,$
$\BI$ being the unit operator acting on the vector spaces.
$r$ acts nontrivially on both vector spaces.    Here
we shall utilize the following $4 \times 4$ matrix realization for $r$:
\Eq{ r={{i\lambda}\over 2}\pmatrix{1 & & &  \cr &-1 & &  \cr &4 &-1& \cr
& & & 1\cr} \quad , \label{rmat}}  $\lambda$ being a real parameter.

Eqs. (\ref{pbptpt}-\ref{pbpgm}) give a one parameter family
of quadratic Poisson structures on ${\cal O}$.
Jacobi identities involving $\tilde p,\;\gamma$ and $\bar \gamma$
are satisfied in part due to $r$ satisfying the classical Yang-Baxter
equations.  (We however found it more convenient to use algebraic
manipulation packages to check them.)
It can also be checked that det($\gamma$) has vanishing
brackets with all observables and hence eqs.
(\ref{pbptpt}-\ref{pbpgm}) are consistent with the
unimodularity condition.  Furthermore, the relations
(\ref{pbptpt}) and (\ref{pbgam})
define skew symmetric brackets, the former being invariant under
hermitean conjugation.

In Sec. 2, we shall show that the quadratic
  algebra on ${\cal O}$ defined by eqs. (\ref{pbptpt}-\ref{pbpgm})
  is preserved under the Poisson action
of the symmetries ${\cal S}$.   There we shall also show that it
is a deformation of the standard
Poincar\'e algebra, the latter being recovered
when $\lambda \rightarrow 0$, and in that limit, (the canonical limit)
 the Lie-Poisson  transformations reduce to canonical transformations.
Like with the Poincar\'e algebra, the algebra described by eqs.
 (\ref{pbptpt}-\ref{pbpgm}) has two Casimir invariants.  One of the
 Casimirs is
the square of the momenta, while the other is the square of a vector
which is a deformation of the Pauli-Lubanski vector.  We can therefore
associate the two Casimirs with ``mass" and ``spin".

In Secs. 3 and 4, we parametrize the symplectic leaves of
${\cal O}$ with variables which we associate with space-time
coordinates $x$,
momenta $\tilde p$ and spin $\gamma_s$.  We give a realization of the algebra
 (\ref{pbptpt}-\ref{pbpgm}) in terms of $x$ and $\tilde p$ in Sec. 3.
 The  Poisson structure for these variables was already
written down in ref. \cite{ssy}.  There it was shown to be preserved
under the Poisson action of ${\cal S}$.  It was also shown to be a
deformation of the canonical symplectic structure for a relativistic
particle, i.e. $\pois{x_\mu}{p_\nu}=\eta_{\mu\nu},\; \eta =$ diag
$(-1,1,1,1)$.    In Sec. 3  (and in the appendix),
we shall write $\gamma$ in terms of space-time
coordinates $x$ and the  momenta $\tilde p$.  This gives an
expression analogous to the orbital angular momentum in the canonical
theory of a relativistic particle.  Only here $\gamma$ is
an infinite series in $x$ and $\tilde p$ (and
$\lambda$), which reduces to the usual expression
for the orbital angular momentum
when $\lambda \rightarrow 0$.
We show that the deformed Pauli-Lubanski vector is zero for this
realization (for any value of $\lambda$), and hence we conclude that we
have a description of a particle with zero spin.
We also remark that if the classical Hamiltonian for the system
is taken to be the momentum squared (i.e, det $\tilde p$)
times a Lagrange multiplier,
then the resulting dynamics is identical to
 that of a free massless
particle (for any value of $\lambda$).  We thus arrive at a one
parameter family of canonically inequivalent descriptions for
the photon.  Upon quantization, the resulting states are expected
to transform under the action of the  quantum Lorentz group.

In Sec. 4, we show how the algebra  (\ref{pbptpt}-\ref{pbpgm})
can be realized when spin is present.  Unlike in the canonical
theory, we find that the spin associated with a particle must have
 nonvanishing Poisson brackets with both the
space-time coordinates and the momenta.
This is a consequence of the result that ${\cal O}$, unlike
${\cal S}$, is not a Poisson-Lie group.  When the mass shell constraint
is taken for the Hamiltonian, the classical spin is
found to have a trivial dynamics, i.e. there is no spin precession,
and this is just as in the canonical theory.\cite{spr}

In Sec. 5, we give a preliminary discussion of the quantization
of the Poisson bracket algebra (\ref{pbptpt}-\ref{pbpgm}),
while concluding remarks are made in Sec. 6.

\section{The Deformed Poincar\'e Algebra}

Here we will examine the two distinct Poisson manifolds
${\cal S}$ and ${\cal O}$, associated respectively with the space
of symmetries and the space of classical observables.  As stated in
the introduction, we shall identify the former
with the six-dimensional Lorentz group
 [or more precisely, its covering group $SL(2,C)$], having
Poisson brackets $\{\;,\;\}_
{\cal S}$ corresponding to that of a Poisson-Lie group
\cite{Tak},\cite{ssy}.
${\cal O}$ will be assumed to be ${\bf R}^4\times SL(2,C)$
with Poisson brackets
$\{\;,\;\}_{\cal O}$ defining a one-parameter deformation of the
Poincar\'e algebra.  The deformed
  algebra is essentially given by eqs. (\ref{pbptpt}-\ref{pbpgm}).

We first review the Poisson structure on ${\cal S}$.
\cite{Tak},\cite{ssy}
For simplicity of notation we shall drop the subscripts on
the Poisson brackets.

\subsection{Symmetries}

Let $g$ be a $2 \times 2$ complex unimodular matrix which we use to
parametrize ${\cal S}$.  For ${\cal S}$ to be a Poisson-Lie group its
Poisson brackets must be compatible with left and
right group multiplication.  Such brackets are:
\Eq { \{\gone,\gtwo\} = [\;r  \;,\;\gone \gtwo\;] \;,\label{lpb} }
where we again utilize tensor product notation, with
$\gone=g\otimes \BI,\;\gtwo=\BI\otimes g $ and the  $r-$ matrix
defined in eq. (\ref{rmat}).  Since the latter is proportional to
$\lambda,$  the group elements have zero
Poisson brackets in the limit $\lambda\rightarrow 0$, and this
once again corresponds to the canonical limit.
The Jacobi identity holds due to the $r-$matrix satisfying
the (modified) classical Yang-Baxter equation.
 (The Leibniz identity for the
Poisson brackets is assumed here and throughout this article.)
It can also be checked from eq.
(\ref{lpb}) that det($g$) has zero Poisson
brackets with all components of $g$ and hence we may consistently
set det$(g)=1$.

The transformations of the observables will involve $g$ as
well as its hermitean conjugate.  Therefore in addition to eq.
(\ref{lpb}) we will also need to know the Poisson brackets of
$g^\dagger$ or
 $\bar g = (g^\dagger)^{-1}$.  For this we demand
that the Poisson structure for $g$ and $ g^\dagger$ is consistent
with complex conjugation, antisymmetry and the Jacobi identity.
All three of these conditions are met  for the following set of relations
\cite{ssy}:
\begin{eqnarray}
 \{\gone,\gbtwo\}&=&[\;r  \;,\;\gone \gbtwo\;] \;,\label{ggb}  \\
 \{\gbone,\gtwo\}&=&[\;r^\dagger\;,\;\gbone \gtwo\;]\;,\label{gbg} \\
 \{\gbone,\gbtwo\}&=&[\;r  \;,\;\gbone \gbtwo\;] \;.\label{gbgb}
\end{eqnarray}
As indicated earlier, the Poisson brackets
(\ref{lpb}) and (\ref{ggb}-\ref{gbgb}) coincide with the classical limit
of the $SL_q(2,C)$ commutation relations given in refs.\cite{AKR},
\cite{Tak}.

We note that eqs. (\ref{lpb}) and (\ref{gbgb}) can be rewritten with
$r$ replaced by $r^\dagger$.  This is  because
$r-r^\dagger$ serves as an adjoint invariant for $SL(2,C)$.  More
specifically,  using the matrix representation (\ref{rmat}) for $r$
we have the identity
\Eq{r-r^\dagger= i\lambda(2\Pi- \BI) \;,  \label{rrdid} }
where $\BI$ is the unit operator (now acting on the entire
tensor product space) and $\Pi$ is the permutation operator,
i.e. $\Pi$ switches the two vector spaces.  Thus for example,
$\gone \Pi=\Pi \gtwo$ and $\gtwo \Pi=\Pi \gone$.

\subsection{Observables}

Here we discuss the Poisson structure on ${\cal O}$ given in eqs.
(\ref{pbptpt}-\ref{pbpgm}) expressed in terms of the ten observables
contained in $\tilde p$ and $\gamma$.  The former transforms as
a Lorentz vector, i.e. it corresponds to the
 $(\frac12,\frac12)$ representation of the Lorentz group,
 and we associate it with momenta, while the latter transforms as the
$(1,0)$ and $(0,1)$ representations and we associate it with  angular
momenta.   We shall here show that the Poisson structure
(\ref{pbptpt}-\ref{pbpgm}) is a deformation
of the standard Poincar\'e algebra and that it is preserved under
the Lie-Poisson action of the Lorentz group.

We first discuss the Poisson structure (\ref{pbptpt})
for the momenta $\tilde p$.  Actually, this structure was already
 given in ref. \cite{ssy}.
 There we wrote $\tilde p$ as a $2 \times 2$
hermitean matrix   \Eq{\tilde
 p=\pmatrix{-p_0+p_3 &  p_1-ip_2\cr p_1+ip_2 & -p_0-p_3} \;,\label{pmat}}
$p_\mu$ being the space-time components.  Under a Lorentz transformation
\Eq{ \tilde p\rightarrow \tilde p'=\bar g\tilde pg^{-1} \;.\label{Ptr} }
The Poisson structure for $\tilde p$ is required to
 be preserved under these transformations, where we assume the
Poisson brackets (\ref{lpb}) and (\ref{ggb}-\ref{gbgb}) for $g$ and $\bar g$.
(As stated earlier, the symmetries $g$ and $\bar g$
are assumed to have zero Poisson brackets with all
observables.)   The Poisson brackets
 are also required to be skewsymmetric, invariant
under hermitean conjugation and satisfy the Jacobi identity.

 A solution to all of the above requirements is eq. (\ref{pbptpt}).
It is easy to check that eq. (\ref{pbptpt}) is preserved under
Lorentz transformations,
\begin{eqnarray}
\pois{\ptone}{\pttwo} \rightarrow  \pois{\ptone'}{\pttwo'}
& &= \pois{\gbone\ptone \gone^{-1}}{\gbtwo\pttwo \gtwo^{-1}}   \cr
& &=   r \ptone'\pttwo' +\ptone'\pttwo'
r^\dagger  -\pttwo' r^\dagger \ptone' -\ptone' r \pttwo'\;.
\end{eqnarray}
Since the $r-$matrix is proportional to $\lambda$,
all of the brackets are zero in the limit $\lambda \rightarrow 0$
and we recover the canonical result.  The skewsymmetry of the bracket
and invariance under hermitean conjugation is also easily checked.

In terms of the space-time
components of $\tilde p$, eq. (\ref{pbptpt}) can be written as
\begin{eqnarray}
\pois{p_i}{p_j}&=&2\lambda \epsilon_{ijk} p_k (p_0+p_3)\;,\quad  \cr
\pois{p_i}{p_0}&=&0\;,\cr & & \qquad\qquad\qquad\qquad i,j,k=1,2,3\;.
\end{eqnarray}
Thus the time component $p_0$ is in the center of the algebra.
Also in the center is the magnitude of
spatial components $\sqrt{p_ip_i}$ and consequently the
invariant mass-squared, i.e. $p^\mu p_\mu =$
det($\tilde p$).  We expect that analogous central
elements appear in the quantum theory, indicating that
simultaneous measurements of the ``energy", magnitude of the ``momentum"
and one spatial component of $p_\mu$ are possible.\cite{OSWZ}

We next take up the Poisson structure of the angular momenta which we
denote by $j$.   $j$ can be represented by a
$2 \times 2$ complex traceless matrix.  Actually,
 we find it more convenient
however to deal with the exponentiation of $j$ which we denote by
$\gamma=e^{i\lambda j}$.  Like $g$, $\gamma$ is unimodular, i.e.
det$(\gamma) = 1$, and hence it is an $SL(2,C)$ matrix.  Our space of
classical observables is thus ${\bf R}^4\times SL(2,C)$.

  Under Lorentz transformations
\Eq{ \gamma \rightarrow \gamma'= g\gamma g^{-1} \;.\label{gamtr} }
The Poisson structure for $\gamma $ is required to
 be covariant under this transformation.
It is also required to be antisymmetric, consistent with the
constraint det$(\gamma) = 1$ and satisfy the Jacobi identity.  A
solution is eq. (\ref{pbgam}).   Under a Lorentz transformation,
\begin{eqnarray}
\pois{\gamone}{\gamtwo} \rightarrow  \pois{\gamone'}{\gamtwo'}
& &= \pois{\gone\gamone \gone^{-1}}{\gtwo\gamtwo \gtwo^{-1}}   \cr
& &= r^\dagger \gamone'\gamtwo' +\gamone'
\gamtwo' r  -\gamtwo'r \gamone'-\gamone'r^\dagger \gamtwo'\;.
\end{eqnarray}
and hence eq. (\ref{pbgam}) is preserved.
{}From eq. (\ref{pbgam})
it can also be checked from that det$(\gamma)$ has zero Poisson
brackets with all components of $\gamma$ and hence we may consistently
set det$(\gamma)=1$.

In addition to eq.
(\ref{pbgam}), we need to specify the Poisson brackets of
$\gamma^\dagger$ or $\bar \gamma = (\gamma^\dagger)^{-1}$.
For this we again demand that it be
  preserved under Lorentz transformations and that it is consistent
with complex conjugation, antisymmetry and the Jacobi identity.
All of these conditions are met  for eq.
(\ref{pbgmgmb}) along with the following relations:
\begin{eqnarray}
 \{\gambone,\gambtwo\}&=&
  r\gambone\gambtwo +\gambone\gambtwo r^\dagger
  -\gambtwo r^\dagger \gambone -\gambone r \gambtwo \;,\label{pbgmbgmb}
\\
 \{\gambone,\gamtwo\}&=&
  r^\dagger \gambone\gamtwo +\gambone\gamtwo r^\dagger
  -\gamtwo r^\dagger \gambone -\gambone r^\dagger \gamtwo \;.
\label{pbgmbgm}
\end{eqnarray}

The remaining Poisson brackets of the observables are between $\tilde p$
and the group variables $\gamma$ and $\bar \gamma$.  For them
we find eq.
(\ref{pbpgm}) along with \Eq{
   \{\ptone,\gambtwo\} =   r\ptone\gambtwo +\ptone\gambtwo r
  -\gambtwo r^\dagger \ptone -\ptone r \gambtwo \;.\label{pbpgmb} }
Using eqs. (\ref{pbgmgmb}) and (\ref{pbpgmb}),
we have checked that $\pois{{\rm det}\;(\gamma)}{\tilde p}=
\pois{{\rm det}\;( \gamma)}{\bar \gamma}= 0$ and therefore that these
Poisson brackets are consistent with the condition of unimodularity.

The Poisson structure of all ten observables is given by eqs.
 (\ref{pbptpt}-\ref{pbpgm})
and (\ref{pbgmbgm}-\ref{pbpgmb}).  [Actually, we only need to specify
eqs. (\ref{pbptpt}-\ref{pbpgm})
as the remaining relations are obtained by hermitean conjugation.]
We have used algebraic manipulation packages to verify the Jacobi
identity for $\tilde p$, $\gamma$ and $\bar \gamma$.

We next show that the algebra generated by $\tilde p,\;\gamma$
and $\bar \gamma$ is a deformation of the
standard Poincar\'e algebra, the latter being recovered in the limit
$\lambda\rightarrow 0$.  For this we substitute
$\gamma=e^{i\lambda j}$ and $\bar \gamma=e^{i\lambda j^\dagger}$
 into the Poisson bracket relations
and expand around $\lambda=0$, keeping only the lowest order
contributions.  As stated earlier, eq. (\ref{pbptpt}) gives
\Eq{\pois{\ptone}{\pttwo}  \rightarrow 0 \;.\label{canl1}}
The lowest order contributions to eqs.
 (\ref{pbgam}) and (\ref{pbgmgmb}) are quadratic in $\lambda$, yielding
\Eq{  \{\jone,\jtwo\}  \rightarrow  2\Pi\;
  (\jtwo-\jone)   \;,\qquad
 \{\jone,\jdtwo\} \rightarrow  0\;,   \label{canl3} }
where we used eq. (\ref{rrdid}).  Lastly, from eq. (\ref{pbpgm}) we get
\Eq{ \{\ptone,\jtwo\} \rightarrow
\ptone \; (2 \Pi - \BI )\;.\label{canl4}}
Eqs. (\ref{canl1}-\ref{canl4}) define the Poincar\'e algebra.  It can
be expressed in a more familiar, i.e.
\begin{eqnarray}
 \{p_\mu,p_\nu\} &=& 0\;,    \\
 \{j_{\mu\nu}, j_{\rho \sigma}\}&=&\eta_{\mu\rho} j_{\nu\sigma}
+\eta_{\nu\sigma} j_{\mu\rho} +\eta_{\mu\sigma} j_{\rho\nu}
+\eta_{\nu\rho} j_{\sigma\mu}  \;,\\
 \{p_{\mu}, j_{\nu \rho}\}&=&\eta_{\mu\rho} p_{\nu}
-\eta_{\mu\nu} p_{\rho} \;,\\
& &\qquad\qquad \eta={\rm diag}(-1,1,1,1)\;,
\end{eqnarray}
 upon applying the matrix representation [cf. eq.
(\ref{pmat})] for $\tilde p$, along with the following representation
for the $2\times2$ complex traceless matrix $j$:
\Eq{j=\pmatrix{-i j_{12} + j_{30} & -ij_{23}-ij_{20}-j_{31}+j_{10}\cr
 -ij_{23}+ij_{20}+j_{31}+j_{10} & i j_{12} - j_{30} }\label{jmat}}

\subsection{Casimirs}

Like the Poincar\'e algebra, the algebra
generated by $\tilde p,\;\gamma$ and $\bar \gamma$ has two central
elements, which we will associate  with ``mass" and ``spin".

With regard to the mass, this classical Casimir function
is identical in form
to that of the Poincar\'e algebra.   (This however isn't the case at
the quantum level \cite{OSWZ}.  To define the latter one normally
introduces a deformed determinant.)  That is,
$p^\mu p_\mu =$det($\tilde p$) is the Casimir function.
{}From eqs.
(\ref{pbptpt}), (\ref{pbpgm}) and (\ref{pbpgmb}), we have that
\Eq{ \pois{{\rm det}\;(\tilde p)}{\tilde p}=
\pois{{\rm det}\;(\tilde p)}{\gamma}=
\pois{{\rm det}\;(\tilde p)}{\bar \gamma}= 0\;.\label{cas1}}
and therefore that it is in the center of the algebra.

With regard to the spin, the second Casimir function
 can be defined as the square
of a new vector $\tilde w$ which we now define:
\Eq{\tilde w={1\over{2\lambda}}\;(
\bar \gamma ^{-1} \tilde p \gamma  - \tilde p)\;.\label{wdef}   }
 $\tilde w$ is a
$2 \times 2$ hermitean matrix, so we can write:
\Eq{\tilde w=\pmatrix{-w_0+w_3 &  w_1-iw_2\cr w_1+iw_2 & -w_0-w_3}\;.}
It is also a deformation of the standard Pauli-Lubanski vector.
 To see this, we substitute
$\gamma=e^{i\lambda j}$ and $\bar \gamma=e^{i\lambda j^\dagger}$
in eq. (\ref{wdef}) and expand around $\lambda=0$, yielding
\Eq{\tilde w = {i\over2} (\tilde p j- j^\dagger \tilde p) +{\cal O}
(\lambda)\;,}  the zeroth order term in $\lambda$ being
the Pauli-Lubanski vector.  Under Lorentz transformations,
$\tilde w$ transforms as $\tilde p$ does, i.e.
$ \tilde w\rightarrow  \tilde w'= \bar g\tilde wg^{-1}$, and
in addition, we find that
its Poisson brackets with the observables $\gamma$ and $\bar \gamma$
are identical in form
to  those of $\tilde p$ with $\gamma$ and $\bar \gamma$, i.e.
\begin{eqnarray}
 \{\wtone,\gamtwo\}&=&  r^\dagger \wtone\gamtwo +\wtone\gamtwo r
  -\gamtwo r^\dagger \wtone -\wtone r^\dagger \gamtwo \;,\label{pbwgm}\\
    \{\wtone,\gambtwo\}&=&  r\wtone\gambtwo +\wtone\gambtwo r
  -\gambtwo r^\dagger \wtone -\wtone r \gambtwo \;.\label{pbwgmb}
\end{eqnarray}
The Poisson brackets $\tilde w$ with $\tilde p$ are given by
\Eq{\pois{\wtone}{\pttwo} =  r\wtone\pttwo +\wtone\pttwo
r^\dagger  -\pttwo r^\dagger \wtone -\wtone r \pttwo \;,\label{pbptwt}}
which in terms of space-time components, eq. (\ref{pbptwt})
 can be expressed as follows:
\begin{eqnarray}
\pois{w_i}{p_j}&=&2\lambda\biggl(
 \epsilon_{ij\ell}  (p_0+p_3) -\delta_{j3}\epsilon_{ik\ell} p_k
\biggr) w_\ell    \;,\cr
\pois{w_i}{p_0}&=&2\lambda \epsilon_{ijk}  p_jw_k \qquad i,j,k=1,2,3\;,
\cr
\pois{w_0}{p_\mu}&=&0 \;.    \label{wpb2}
\end{eqnarray}
We then find that, in analogy to eq. (\ref{cas1}),
\Eq{ \pois{{\rm det}\;(\tilde w)}{\tilde p}=
\pois{{\rm det}\;(\tilde w)}{\gamma}=
\pois{{\rm det}\;(\tilde w)}{\bar \gamma}= 0\;,\label{cas2}}
and hence that $w_\mu w^\mu =$ det($\tilde w$) is a classical
 Casimir function.

We expect that there are no additional independent Casimir functions
and therefore that the symplectic leaves in ${\cal O}$ are
eight-dimensional, just as is the case with the Poincar\'e algebra.  In
the two sections which follow, we shall show how to parametrize the
symplectic leaves with variables which one can naturally associate
with position, momenta and spin.

For completeness we compute the Poisson brackets for $\tilde w$ with
itself.  From eqs. (\ref{pbwgm}-\ref{pbptwt}), we find
\Eq{\pois{\wtone}{\wttwo} =  r\wtone\wttwo +\wtone\wttwo
r^\dagger  -\wttwo r^\dagger \wtone -\wtone r \wttwo
-i\Pi (\wtone\pttwo-\wttwo \ptone) \;,        \label{pbwtwt}}
or in terms of the space-time components of $\tilde w$,
\begin{eqnarray}
\pois{w_i}{w_j}&=&\epsilon_{ijk}\biggl(p_0w_k -w_0p_k +2\lambda
  (p_0+p_3)p_k \biggr)    \;,\cr
\pois{w_0}{w_i}&=& \epsilon_{ijk}  p_jw_k \qquad i,j,k=1,2,3
\;.    \label{wwb2}
\end{eqnarray}
{}From eqs. (\ref{wpb2}) and (\ref{wwb2}), we deduce
that the set of commuting operators
in the quantum theory can be enlarged to those associated with
$$p_0,\;p_3,\;w_0\; {\rm and} \;w_ip_i\;,$$ in addition to the
two Casimirs $p_\mu p^\mu$ and  $w_\mu w^\mu$.

[We note that in ref. \cite{OSWZ}, a set of commuting
operators for the spinless particle
was found for a similar system.  The set contained operators
associated with  $p_0$, $p_3$ and the third component
of angular momentum.  All of these variables were shown to
have a discrete quantum spectrum for the case of a particle with nonzero mass.
(From ref. \cite{ssy}, we surmise that such a particle is not
free, but instead has a nontrivial interaction with the space-time.)
With regard to the variables $p_0$ and $p_3$,
we expect that a similar spectrum will occur
for us.    We do not know what the
third component of angular momentum corresponds to in our formalism,
nor do we know if $w_0\; {\rm and} \;w_ip_i\;$ can be included in the
set of commuting operators of ref. \cite{OSWZ}.]

\section{Spin Zero Realization}

Here we discuss a realization of the deformed Poincar\'e
algebra (\ref{pbptpt}-\ref{pbpgm}) in terms of space-time
coordinates $x_\mu$ and the momenta $p_\mu$.  The realization,
is based on the system described in \cite{ssy}, where we deformed
the canonical symplectic structure for a relativistic particle.  As we
shall see, this realization has
 the Casimir  det$(\tilde w)  $ equal to zero, and we therefore can
associate it with the description of a (deformed)
 spinless relativistic particle.   Actually, here
we get the even stronger constraint that $\tilde w=0$, or
from  (\ref{wdef}),
\Eq{ \tilde p \gamma  =\bar \gamma \tilde p\;.\label{weq0}}
This is analogous to what is obtained
in the {\it canonical theory} of spinless particles, where
all of the components of the Pauli-Lubanski vector vanish.

The momentum matrix $\tilde p$ was introduced previously in eq.
(\ref{pmat}).   With regard to the space-time coordinates $x_\mu$,
 we find it convenient to define the $2 \times 2$ hermitean
matrix   \Eq{x=\pmatrix{-x_0-x_3 & -x_1+ix_2\cr -x_1-ix_2 & -x_0+x_3}\;
 \;.\label{xmat}}
In contrast to $\tilde p$, $x$ transforms according to
\Eq{x\rightarrow x'= gx\bar g^{-1} \;.\label{xtr} }

As stated previously, the Poisson structure for $x$ and $\tilde p$ is
required to be
a deformation of the canonical Poisson brackets for a relativistic
particle.   In addition it should be preserved under the Lie-Poisson
action of the Lorentz group,
satisfy the Jacobi identity, and be hermitean.
 The following Poisson brackets are
 consistent with all of the above conditions:
\begin{eqnarray}
  \pois{\xone}{\xtwo}&=& r \xone\xtwo +\xone\xtwo
r^\dagger  -\xtwo r \xone -\xone r^\dagger \xtwo \;,\label{pbxx}\\
 \pois{\xone}{\pttwo}&=& r \xone\pttwo +\xone\pttwo
r^\dagger  -\pttwo r \xone -\xone r^\dagger \pttwo
 - {\Pi}(\fone^\dagger +\ftwo) \;,\label{pbxpt}
\end{eqnarray}
along with eq. (\ref{pbptpt}).    $f$ is a $2\times 2$ complex matrix.
For eq. (\ref{pbxpt}) to be preserved under Lorentz transformations,
it must transform like $\gamma$, i.e.
\Eq{  f \rightarrow  g      f g^{-1} \;.\label{flptr} }
$f$ must be a function of $\lambda$.  This is since
in order to recover the canonical relations, we need that $f$
tends to the unit matrix $\BI$ when $ \lambda\rightarrow 0,$
while it cannot be $\BI$ for all $\lambda$ because some work shows
that it would not then satisfy the Jacobi identity.  In this regard, the
 issue of the Jacobi identity was only partially
 addressed in \cite{ssy}.  Here we find (with
the aid of algebraic manipulation packages) that the Jacobi
identity involving the position and momentum variables
holds provided that $f$ satisfies the following Poisson
brackets with $x$ and $\tilde p$:
\begin{eqnarray}
 \{\ptone,\ftwo\}&=&   r^\dagger \ptone\ftwo +\ptone\ftwo r
  -\ftwo r^\dagger \ptone -\ptone r^\dagger \ftwo +i\lambda\ptone
\ftwo    \;,\label{pbpf}
  \\   \{\xone,\ftwo\}&=&  r^\dagger \xone\ftwo +\xone\ftwo r^\dagger
  -\ftwo r \xone -\xone r^\dagger \ftwo -i\lambda\xone
\ftwo    \;.\label{pbxf}
\end{eqnarray}

In the above, it appears that we have enlarged the phase space
spanned by $x$ and $\tilde p$ to also include the variables $f$.
However it is not necessary to regard $f$ as independent variables.
Rather,  it is possible to express $f$ in terms of $x$ and $\tilde p$,
while still being consistent with Poisson brackets
 (\ref{pbpf}) and (\ref{pbxf}), as well as with the canonical limit
$f\rightarrow \BI$.
This was done in ref. \cite{ssy}, where we wrote $f$ according to:
 \Eq{f= \exp \{i\lambda
J \}\;,\quad {\rm where} \quad \sin \lambda J=\lambda x\tilde p\;.
\label{feiJ}   }
By this we meant the following Taylor expansion:
\begin{eqnarray}
f&=&\BI + i\lambda x\tilde p -{1\over 2} (\lambda x\tilde p)^2
 - {1\over 8} (\lambda x \tilde p)^4 \; ... \cr
& &\cr
 & &\quad...\; - \frac{1\cdot 3\cdot 5\cdot\;...\;\cdot (2n-3)}
{2\cdot 4\cdot 6\cdot\;...\;\cdot 2n}
(\lambda x\tilde p)^{2n}\;...  \label{polf}
\end{eqnarray}
It is easily seen that this expression transforms as in eq.
(\ref{flptr}) and that it
tends to the unit matrix $\BI$ when $ \lambda\rightarrow 0.$
We show that it agrees with eq. (\ref{pbpf}) in the appendix.
In a similar manner, it can be shown to agree with eq. (\ref{pbxf}).
Also in a similar manner, we can use eq. (\ref{feiJ}) to
obtain the brackets for $f$ with itself and with $f^\dagger$.  We find:
\begin{eqnarray}
 \{\fone,\ftwo\}&=&  r^\dagger \fone\ftwo +\fone\ftwo r
  -\ftwo r \fone -\fone r^\dagger \ftwo \;,\label{pbff}  \\
 \{\fone,\fdtwo\}&=&   r \fone\fdtwo +\fone\fdtwo r
  -\fdtwo r \fone -\fone r \fdtwo \;.  \label{pbffd}
\end{eqnarray}
We have checked that these relations are consistent with the Jacobi
identity.  We note that the former Poisson bracket
is identical in form to the Poisson brackets
(\ref{pbgam}) for $\gamma$, the only difference between them being that while
$\gamma$ appearing in eq. (\ref{pbgam})
is unimodular, $f$ is not.  In this regard, we note that we
are not allowed
to set det$(f)$ equal to one, because although
 det$(f)$ has zero brackets
with $f$ and $f^\dagger$ [which follows from eqs. (\ref{pbff})
and (\ref{pbffd})], it does not have zero brackets with the coordinates
or the momenta.  Instead from eq. (\ref{pbpf}),  we find
\begin{eqnarray}
\pois{x}{{\rm det}\;(f)}&=&-2i\lambda x\; {\rm det}\;(f) \;,\label{xf}\cr
\pois{\tilde p}{{\rm det}\;(f)}&=&2i\lambda \tilde p\; {\rm det}\;(f) \;.
\end{eqnarray}

Alternatively, we can easily define a unimodular matrix from $f$
by simply dividing by $\sqrt{{\rm det}( f)}\;$.  This is in fact how we
make the identification with $\gamma$.  Thus we set
\Eq{\gamma = \frac{f} {\sqrt{{\rm det}( f)}} \;.\label{gfid} }
It then follows that if we again write $\gamma=e^{i\lambda j}$
and use eq. (\ref{feiJ}),  then  $j$ corresponds to the
traceless part of $J$,   \Eq{j = J-\frac 1 2 \BI \;{\rm Tr }\;J\;. }
Upon keeping terms linear in $\lambda$ in the expression (\ref{polf}) for
we see that $j$ reduces to the usual expression for orbital angular
momentum in the canonical limit $\lambda\rightarrow 0$.

Using eq. (\ref{gfid}) we can recover the Poisson brackets
(\ref{pbgam}-\ref{pbpgm}) starting from the brackets for $f$
[specifically,
(\ref{pbpf}), (\ref{pbff}) and (\ref{pbffd})].  Since we also
have eq. (\ref{pbptpt}), we obtain a realization of the entire
deformed Poincar\'e algebra.  Using eqs.
(\ref{pbxf}) and (\ref{xf}), we can further obtain
 the brackets between the space-time coordinates $x$ and $\gamma$: \Eq{
    \{\xone,\gamtwo\} =   r^\dagger \xone\gamtwo +\xone\gamtwo r^\dagger
  -\gamtwo r \xone -\xone r^\dagger \gamtwo \;.}

It now only remains to show that
 the condition (\ref{weq0}) is satisfied implying that the
deformed Pauli-Lubanski vector is zero .   This follows
using results from the appendix, specifically eqs.
(\ref{pssp}) and (\ref{pccp}), which lead to
\Eq{\tilde p e^{i\lambda J} =  e^{i\lambda J^\dagger}  \tilde p\;,}
If we now divide both sides of this equation by
$\sqrt{{\rm det}( f)}\;$, we get eq. (\ref{weq0}) and hence the
 description of a (deformed) spinless particle.  We therefore
 associate the expression for $\gamma$ in terms of $x$ and $\tilde p$
with the `orbital angular momentum'.

In \cite{ssy}, we made some remarks concerning the dynamics for such
a system.  There we showed that if the standard mass shell constraint
is chosen for the Hamiltonian function, i.e.
\Eq{ H=\alpha( \det(\tilde p) - m^2) \;,\label{ham}}
then it, along with the Poisson brackets
(\ref{pbptpt}), (\ref{pbxx}) and (\ref{pbxpt}), yields a nontrivial
interaction of the  particle with the space-time when $\lambda$ and
$m$ are different from zero.  ($\alpha$
denotes a Lagrange multiplier).   In \cite{ssy} we solved for the particle
trajectory and found that it originates and terminates at
singularities.  The particle has a lifetime equal to
\Eq{\bigg| \frac{{\rm Tr}(\gamma  p)}
{\lambda \det (\tilde p)}\bigg|\;,\label{lt}}
 where $p=\sigma_2 \tilde p^T \sigma_2\;,$
 $T$ denoting transpose and $\sigma_2$ being the second Pauli matrix.
Thus when $\lambda$ and
$m$ are different from zero, it appears that the dynamics corresponds
to  a kind of `virtual' particle which is off its `physical' mass shell.
(By `physical' mass shell we mean that the square of the {\it canonical}
momentum - not $p_\mu$ - is equal to $m^2$.)

We note that the expression (\ref{lt}) is singular for the case of a
massless particle, i.e. $m=0$, indicating that such a particle has
an infinite lifetime.  Actually, we
can show that the Hamiltonian (\ref{ham}) with $m=0$ [along
with the Poisson brackets
(\ref{pbptpt}), (\ref{pbxx}) and (\ref{pbxpt})]
 describes a free photon (or any massless particle)
for  arbitrary values of $\lambda$.
For this we compute  the Hamilton equations of motion:
\begin{eqnarray}
\dot x ={d\over {d\tau}} x &=& \alpha
 \{x, \det (\tilde p) \}=-\alpha (f p + p f^\dagger)  \;,  \\
\dot {\tilde p} ={d\over {d\tau}} \tilde p &=& \alpha
 \{\tilde p, \det (\tilde p) \}=0     \;. \end{eqnarray}
It then follows that for $m=0$,
\Eq{\dot x \tilde p = -\alpha (f p + p f^\dagger) \tilde p  =0  \;,
\label{xdtp} }
where we have used the mass shell constraint $\tilde p p=p \tilde p =$
${\rm det}(\tilde p)  \BI=0$, in addition to the infinite series
 expression
 eq. (\ref{polf}) for $f$.    The traceless part of eq. (\ref{xdtp})
is equivalent to $\dot x_\mu p_\nu
-\dot x_\nu p_\mu=0$ and hence
\Eq{\dot x_\mu = \kappa p_\nu \;.\label{frph}}
The trace of eq. (\ref{xdtp}) implies that
$\dot x^\mu p_\mu =0$, and thus gives no constraint
on the proportionality constant $\kappa$.  We therefore arrive at a
{\it  free} light-like trajectory.  Furthermore,
since $\lambda$ is arbitrary,
we get an entire family of canonically inequivalent Hamiltonian
descriptions of a photon trajectory.
 Upon quantization, the resulting states are expected
to transform covariantly under the action of the quantum Lorentz group.

Actually, to truly
 describe a photon, we should introduce a spin and check that
its equation of motion is the usual one.   That is, there should be
no classical spin precession.  Spin will be introduced
in the next section.  There we indeed find that there is no precession
of the classical  spin (even for the case $m \ne 0$).

\section{Inclusion of Spin}

In the canonical theory, spin is introduced as an additional term
in the angular momentum.  This term is defined to have
 zero Poisson brackets
with coordinates and momenta, and it is the only term in the angular
momentum which contributes to the Pauli-Lubanski vector.
We shall look for an analogous prescription for including spin in
 our deformed Poincar\'e algebra.  It should reduce to the canonical
 prescription  in the limit $\lambda \rightarrow 0$.  In this regard,
 there is of  course no unique
prescription for including spin in the deformed Poincar\'e algebra.
In what follows,
our choice  shall be to  multiply  the $SL(2,C)$ matrix
$\gamma$ obtained in the previous
section on the right by another $SL(2,C)$ matrix $\gamma_s$.  The
former is to be regarded as the orbital angular momentum,
  while $\gamma_s$ plays the role of `spin'.   Thus
\Eq{ \gamma \rightarrow  \gamma \gamma_s
 = \frac{f\gamma_s} {\sqrt{{\rm det}( f)}}\;.\label{spin}}
To get back the canonical prescription, i.e. $j\rightarrow j+s$,
 when $\lambda\rightarrow 0$, we can
take $\gamma_s=e^{i\lambda s}$, $s$ being a traceless complex matrix.
Furthermore, using eq. (\ref{spin}) we get that the `spin' $\gamma_s$,
and not the `orbital angular momentum'
$\gamma$  [subject to eq. (\ref{weq0})],
contributes to the deformed Pauli-Lubanski vector eq.
(\ref{wdef}), analogous to what happens in the canonical theory.
 [This would not have been the case if instead of eq.
(\ref{spin}), we had  multiplied $\gamma$ on the left by $\gamma_s$.]
We next show that {\it unlike} in the canonical theory, the spin has nonzero
Poisson brackets with momentum and position, and that this is a
consequence of the fact that the space spanned by the matrices
$\gamma$ does not form a Poisson-Lie group (unlike the space spanned
by matrices $g$).

Because ${\cal S}$ is a Poisson-Lie group the Poisson structure for the
$SL(2,C)$ matrices $g$ is preserved under left or right
group multiplication, i.e. the Poisson brackets (\ref{lpb}) are compatible
with the group product\cite{D}-\cite{AM}.
As we show below, the analogous statement does not, however, apply
to the Poisson structure
for the $SL(2,C)$ matrices $\gamma$.
That is, the Poisson brackets (\ref{pbgam})
are not compatible with group multiplication and hence the space
spanned by $\gamma$ is not a Poisson-Lie group.

To see that the Poisson structure for
 symmetries is preserved under group multiplication,
 one defines a variable
$g'\in SL(2,C)$ which  satisfies the same relations as $g$,
\Eq { \{\gone',\gtwo'\} = [\;r  \;,\;\gone' \gtwo'\;] \;, }
 in addition to $ \{\gone',\gtwo\} = 0$.  Then under
right multiplication $g \rightarrow g g'$, we get
\begin{eqnarray}
 \{\gone, \gtwo\}\rightarrow
 \{\gone \gone',\gtwo \gtwo'\}& = &[\;r  \;,\;\gone\gtwo\;]
 \gone' \gtwo' +   \gone\gtwo[\;r  \;,\;\gone' \gtwo'\;] \cr
& = &[\;r  \;,\;(\gone\gone')( \gtwo\gtwo')\;]  \;,
\end{eqnarray}
and hence that the Poisson structure given by (\ref{lpb})
is preserved.

Let us now try the same thing for the observables $\gamma$.
We define $\gamma
_s\in SL(2,C)$ to satisfy the same relations as $\gamma$, \Eq{
\pois{\gamone_s}{\gamtwo_s}= r^\dagger \gamone_s\gamtwo_s +\gamone_s\gamtwo_s
r  -\gamtwo_s r \gamone_s -\gamone_s r^\dagger \gamtwo_s \;,
\label{pbgams}}
in addition to $ \{\gamone_s,\gamtwo\} = 0$.  But it is easily
checked that this Poisson structure is not preserved under
right  (or left) multiplication $\gamma \rightarrow
\gamma \gamma_s$, and therefore that
 $ f\gamma_s/\sqrt{{\rm det}(f)}$ does not give a realization
of the relations (\ref{pbgam}).

To proceed further we shall drop the assumption that
the product space $\{\gamma\}\times\{\gamma_s\}$ has a product
Poisson structure, i.e. we drop the assumption that
the spin has zero Poisson brackets with the
orbital angular momentum $\gamma$
(and hence also with the coordinates and
momenta) in the deformed theory,
i.e.  $ \{\gamone_s,\gamtwo\}\ne 0$.  Instead we take \Eq{
\pois{\gamone_s}{\gamtwo}= r^\dagger \gamone_s\gamtwo +\gamone_s\gamtwo
r^\dagger -\gamtwo r^\dagger \gamone_s -\gamone_s r^\dagger \gamtwo \;,
\label{pbgspin}}  or equivalently,  \Eq{
\pois{\gamone}{\gamtwo_s}= r \gamone\gamtwo_s +\gamone\gamtwo_s
r -\gamtwo_s r \gamone -\gamone r \gamtwo_s \;.}
It then can be checked
 that the product $\gamma \gamma_s$ carries the same
Poisson structure as $\gamma$ [i.e. eq. (\ref{pbgam})], \Eq{
\pois{\gamone \gamone_s}{\gamtwo \gamtwo_s}=
 r^\dagger (\gamone\gamone_s)(\gamtwo\gamtwo_s) + (\gamone
 \gamone_s)(\gamtwo\gamtwo_s )r  -
(\gamtwo\gamtwo_s) r(\gamone \gamone_s)
-(\gamone\gamone_s) r^\dagger(\gamtwo \gamtwo_s) \;. }
Thus, in this sense we can preserve the Poisson structure for the
observables.   We note that the Poisson bracket relations
 for $\gamma$ and $\gamma_s$ are identical in form to those of
 $\gamma$ and $\gamma^\dagger$.  Then since the Jacobi identity holds
for the latter variables, it must also hold for $\gamma$ and $\gamma_s$.
In addition, we have that $\pois{{\rm det}\;(\gamma)}{\gamma_s}=
\pois{{\rm det}\;(\gamma_s)}{ \gamma}= 0$ and therefore
the Poisson brackets (\ref{pbgspin})
are consistent with the unimodularity of
both $\gamma$ and $\gamma_s$.  We note that eq.
(\ref{pbgspin}) is also consistent with the canonical theory,
because if we write $\gamma=e^{i\lambda j}$
and $\gamma_s=e^{i\lambda s}$ then to lowest order in $\lambda$,
we get that $s$ has zero Poisson brackets with $j$.

As stated before, in the canonical theory the spin has zero Poisson
brackets with the momenta.
Here, however, if the Poisson bracket (\ref{pbpgm}) is
to be preserved under eq. (\ref{spin}), we need that
$ \{\gamone_s,\pttwo\}\ne 0$.  Specifically, we need
\Eq{\pois{\ptone}{\gamtwo_s}= r^\dagger \ptone\gamtwo_s +\ptone\gamtwo_s
r -\gamtwo_s r^\dagger \ptone -\ptone r \gamtwo_s \;.\label{pbspp}}
For then   \Eq{   \{\ptone,\gamtwo\gamtwo_s\}=
  r^\dagger \ptone(\gamtwo\gamtwo_s) +\ptone(\gamtwo\gamtwo_s) r
  -(\gamtwo\gamtwo_s) r^\dagger \ptone -
  \ptone r^\dagger (\gamtwo\gamtwo_s )\;,}  and the relation
 (\ref{pbpgm}) is preserved.  We have checked that eq.
(\ref{pbspp}) is consistent with the Jacobi identity for
$\gamma$, $\tilde p$ and $\gamma_s$.  We also verified that
$\pois{{\rm det}\;(\gamma_s)}{ \tilde p}= 0$ and that $\tilde p$
has zero brackets with $s$ in the limit $\lambda \rightarrow 0$.

In the above [cf. eq. (\ref{pbspp})], we found that the spin $\gamma_s$
does not have zero Poisson brackets with
the momenta.  It also doesn't have zero Poisson brackets with
the position.
That $ \{\gamone_s,\xtwo\}\ne 0$ is easily seen because if
it were not so we would not then be able to recover the correct brackets
(\ref{pbgspin}) for $\gamma_s$ with
the orbital angular momenta $\gamma$ [given as  a function
of $x\tilde p$ in eqs. (\ref{polf}) and (\ref{gfid})].
What works instead is
\Eq{\pois{\xone}{\gamtwo_s}= r \xone\gamtwo_s +\xone\gamtwo_s r^\dagger
 -\gamtwo_s r \xone -\xone r^\dagger \gamtwo_s \;.\label{pbspx}}
{}From it and eq. (\ref{pbspp}) we find that
 \Eq{ \pois{(\xone\ptone)^n}{\gamtwo_s}= r (\xone\ptone)^n\gamtwo_s
+(\xone\ptone)^n\gamtwo_s
r -\gamtwo_s r (\xone\ptone)^n -(\xone\ptone)^n r \gamtwo_s \;.}
Hence the Poisson brackets between $\gamma_s$ and
{\it any}  polynomial function  of $x\tilde p$ has precisely the
same form as the  Poisson brackets between
$\gamma_s$ and $\gamma$.  Thus
 \Eq{ \pois{\fone}{\gamtwo_s}= r \fone\gamtwo_s
+\fone\gamtwo_s
r -\gamtwo_s r \fone -\fone r \gamtwo_s \;.\label{pbfgs}}
  From previous arguments we then  also
know that  $\pois{{\rm det}\;(f)}{\gamma_s}=0$.
The Poisson brackets (\ref{pbgspin}) between $\gamma$ and $\gamma_s$
are recovered by dividing both sides of
eq. (\ref{pbfgs}) by $\sqrt{{\rm det}( f)}\;.$
{}From eq. (\ref{pbspx}) it also follows that $\pois{\xone}{\stwo}
\rightarrow 0$ when $\lambda \rightarrow 0$ and hence
we recover the usual canonical limit.

To summarize, we have obtained a realization of the deformed
Poincare algebra defined by eqs. (\ref{pbptpt}-\ref{pbpgm})
with det $(\tilde w)\ne 0$.
(Actually for this we also need to have the Poisson brackets
between $\gamma_s$ and $\bar \gamma$.
 We shall assume that a consistent set of such brackets exist.)
Thus unlike the previous section, we now have a particle
with spin.  We get back the canonical description of a spinning
particle when $\lambda \rightarrow 0$.  Under Lorentz transformations,
$\gamma_s$ must transform as does $\gamma$, i.e. $\gamma_s \rightarrow
g\gamma_s g^{-1}$.  It is easy to check that all Poisson brackets
with $\gamma_s$ are preserved under such transformations and once
again that the Lorentz group induces a Lie-Poisson action on the
space of observables.  (Here we assume as usual that the classical
observables have zero Poisson brackets with the classical
symmetries.)

Finally, we remark about the spin dynamics.  For this purpose
 it is of interest to compute the Poisson
bracket of $\gamma_s$ with det$(\tilde p)$.  From eq.
(\ref{pbspp}), we find that this Poisson bracket
 vanishes.  This means that if
the Hamiltonian function for the system is once again chosen to
be the usual mass shell constraint (\ref{ham}), then the classical
spin has a trivial dynamics (i.e., there is no precession).
This is just as in the canonical formulation of a classical spinning
particle.\cite{spr}
In other words, the particle interaction which is present due to the
highly nontrivial Poisson structure does not affect the spin.

\section{Towards Quantization}

Here we make some preliminary remarks concerning quantization.
We plan to address this issue more fully in a subsequent publication.

There exists a standard quantization scheme (deformation
quantization) which can be applied for
 the symmetries which takes ${\cal S} $ to a Hopf algebra,
specifically, $SL_q(2,C)$.\cite{Tak}  We remark on this first.
We then comment on a possible quantization of the
classical observables.   The system which results appears to be different
from $q-$Poincar\'e algebras discussed previously in the literature.
[3-7]

With regard to
the symmetries, one standardly replaces $g$ by an
 $SL_q(2,C)$ matrix which we denote by  $T$ and
the Poisson brackets (\ref{lpb}) by the corresponding quantum
commutation
relations.  The matrix elements in $T$ are constrained by the
condition that its `deformed' determinant is equal to one, and this
is the analogue of the unimodularity condition on $g$.   The
commutation relations are given in terms of a quantum $R$ matrix,
satisfying the usual quantum Yang-Baxter equations, and can be
written according to
\Eq{\Rontw \Tone \Ttwo = \Ttwo \Tone \Rontw \;.\label{RTT1}}
This algebra can presumably be realized on the space ${\cal S}$ of
classical symmetries with the use of a star product.
In order to recover the correct classical limit one only needs that
$\Rontw\rightarrow \BI - i\hbar r +{\cal O}(\hbar^2)$ when $\hbar
\rightarrow 0$.  In addition to eq. (\ref{RTT1}), one needs the quantum
analogues of the Poisson brackets
(\ref{ggb}) and (\ref{gbgb}).  For this we introduce
another $SL_q(2,C)$ matrix $\bar T$ which in analogy to the classical
observable $\bar g$
is defined by $\bar T ={T^\dagger}^{-1}$.  Then along with
eq. (\ref{RTT1}), we write
\Eq{\Rontw \Tone \Tbtwo = \Tbtwo \Tone \Rontw \;,\label{RTT2}}
\Eq{\Rontw \Tbone \Tbtwo = \Tbtwo \Tbone \Rontw \;,\label{RTT3}}
which reduces to eqs.
(\ref{ggb}) and (\ref{gbgb}) when $\hbar \rightarrow 0$.
By switching vector space indices $1$ and $2$, we see that we
can replace $\Rontw$ in eqs.
(\ref{RTT1}) and (\ref{RTT3}) by $\Rtwon^{-1}$.
Thus $\Rontw\Rtwon$ must
commute with $\Ttwo\Tone$ and $\Tbtwo\Tbone$.  This is analogous to
the statement that
$r-r^\dagger$ is an adjoint invariant in the classical theory.

Concerning the quantum observables, we shall require that its algebra is
preserved under the action of the quantum symmetries, in analogy to
what happens in the classical theory.  We also want that this algebra
is consistent with the classical Poisson bracket algebra in the limit
$\hbar \rightarrow 0$.  A quantum algebra for the momenta
was already given in ref. \cite{AKR} which is consistent with these
properties, so we will adopt it here.  There
one replaces the classical variable
$\tilde p$ by a $2\times 2$ matrix $P$ whose elements are
operator-valued.  The Poisson brackets (\ref{pbptpt}) are replaced
by what were refered to as reflection equations,
\Eq{\Rontw \Pone \Rontw^{-1} \Ptwo = \Ptwo \Rtwon^{-1} \Pone \Rtwon \;.
\label{RPRP}}
With this choice, one can easily obtain the correct classical limit.
For this one notes that
using the matrix expression (\ref{rmat}) for $r$, one gets that
$\Rtwon\rightarrow \BI + i\hbar r^\dagger +{\cal O}(\hbar^2)$ when $\hbar
\rightarrow 0$.  Furthermore, as desired, the commutation relations
(\ref{RPRP}) are preserved under $SL_q(2,C)$ transformations.
Here in analogy to eq. (\ref{Ptr}), one assumes that $P$ transforms
as a vector under the quantum Lorentz group, i.e.
\Eq{ P\rightarrow P'=\bar T P T^{-1} \;, \label{capptr}}
and that the matrix elements of $T$ and $\bar T$
commute with those of $P$.
Then using the relations (\ref{RTT1}-\ref{RTT3}), one gets that the
left hand side of (\ref{RPRP}) transforms according to
\begin{eqnarray}
\Rontw \Pone \Rontw^{-1} \Ptwo \rightarrow
\Rontw \Pone' \Rontw^{-1} \Ptwo'& = &
\Tbtwo \Tbone \Rontw \Pone \Rontw^{-1} \Ptwo \Tone^{-1} \Ttwo^{-1} \cr
&=& \Tbtwo \Tbone \Ptwo \Rtwon^{-1} \Pone \Rtwon \Tone^{-1} \Ttwo^{-1} \cr
&=&  \Ptwo' \Rtwon^{-1} \Pone' \Rtwon  \;, \label{copr}
\end{eqnarray}
and hence that eq. (\ref{RPRP}) is preserved.
Here we may assume that the quantum matrix $P$ is hermitean,
analogous to the fact that the classical matrix
$\tilde p$ is hermitean.  It is easy to check that this is consistent
with the transformation property (\ref{capptr}).  It is also consistent
with the commutation relation (\ref{RPRP}) provided that we have the
following condition on the quantum $R-$matrix.
\Eq{\Rontw^\dagger=\Rtwon\;.   \label{Rdagg}}

It remains to specify the quantum analogues of Poisson brackets
(\ref{pbgam}-\ref{pbpgm}).  For this we associate the classical
observables $\gamma$ and $\bar \gamma$ with operator-valued $2\times 2$
matrices $\Gamma$ and $\bar \Gamma$, which in analogy to
(\ref{gamtr}) transform  as
\Eq{ \Gamma\rightarrow \Gamma'= T \Gamma T^{-1} \;,\quad
\bar \Gamma\rightarrow \bar \Gamma'=\bar T \bar\Gamma \bar T^{-1} \;.
 } under the action of $SL_q(2,C)$.  As with $t$,
we can assume that its matrix elements in $T$ are constrained by the
condition that its `deformed' determinant is equal to one, in analogy
to the unimodularity condition on $\gamma$.  We propose that the
$\Gamma$'s satisfy the following commutation relations amongst themselves
 and with $P$:
\begin{eqnarray}
\Rtwon^{-1} \Gone \Rtwon \Gtwo
&=&\Gtwo \Rontw \Gone \Rontw^{-1} \;,    \label{RGRG} \\
\Rontw \Gone \Rontw^{-1} \Gbtwo
&=&\Gbtwo \Rontw \Gone \Rontw^{-1} \;,             \\
\Rtwon^{-1} \Pone \Rtwon \Gtwo
&=&\Gtwo \Rtwon^{-1} \Pone \Rontw^{-1} \;. \label{RPRG}
\end{eqnarray}
It can be checked that from these relations
one recovers the correct quadratic algebra, i.e. eqs.
(\ref{pbgam}-\ref{pbpgm}), as $\hbar \rightarrow 0$.
Also, using eqs. (\ref{RTT1}-\ref{RTT3}) and the assumption
that the matrix elements of $T$ and $\bar T$ commute with those of
$\Gamma$ and $\bar \Gamma$, it can be checked that
the commutation relations (\ref{RGRG}-
\ref{RPRG}) are preserved under $SL_q(2,C)$ transformations.
The procedure is analogous to that used in eq. (\ref{copr}).

If we define $\bar \Gamma ={\Gamma^\dagger}^{-1}$ in analogy to what
was done in the
classical theory, then by taking the hermitean conjugates of eqs.
(\ref{RGRG}-\ref{RPRG}) we get the quantum analogues of the
Poisson brackets (\ref{pbgmbgmb}-\ref{pbpgmb}).  Using (\ref{Rdagg}),
we get
\begin{eqnarray}
\Rontw \Gbone \Rontw^{-1} \Gbtwo
&=&\Gbtwo \Rtwon^{-1} \Gbone \Rtwon \;,  \\
\Rtwon^{-1} \Gbone \Rtwon \Gtwo
&=&\Gtwo \Rtwon^{-1} \Gbone \Rtwon \;,             \\
\Rontw \Pone \Rontw^{-1}\Gbtwo
&=&\Gbtwo \Rtwon^{-1} \Pone \Rontw^{-1} \;,
\end{eqnarray}
which correspond to eqs.  (\ref{pbgmbgmb}-\ref{pbpgmb}) when $\lambda
\rightarrow 0$.

Although the set of symmetry operators $\{T\}$ defines a Hopf
algebra, the same does not seem to be the case for the set of quantum
operators $\{\Gamma\}$.  In this regard we do not know how to define
a coproduct for the latter.  This is not too surprising since the set
of classical variables $\{\gamma\}$ did not define a Poisson-Lie
group.   For this reason it appears that our quantum algebra defined
in eqs. (\ref{RPRP}) and (\ref{RGRG}-\ref{RPRG}) differs from those
given previously.[3-7]   Eqs. (\ref{RPRP}) and
(\ref{RGRG}-\ref{RPRG}) define a set the quadratic  commutation
relations between the quantum mechanical
observables $P$,$\Gamma$ and $\bar \Gamma$, which nevertheless are
preserved under the action of a Hopf algebra.
Whether or not  it is necessary to impose
 higher order relations remains to be checked.
{}From the observables it is possible to construct the
quantum mechanical deformed Pauli-Lubanski vector, and also
Casimir operators corresponding to mass and spin.  It should then be
possible to look for eigenvectors of these operators along with the
quantum analogues of $p_0,\;p_3,\;w_0\; {\rm and} \;w_ip_i\;.$
We can also hope to obtain realizations of the quantum algebra
for the cases of a spinless and spinning relativistic particle
in a manner similar to what was done in Secs. 3 and 4.

\section{Concluding Remarks}

Here we outline additional future avenues of research.

We have obtained a deformation of the Poincar\'e algebra
which is covariant with respect to the Lie-Poisson  action of the
 Lorentz group.   It is of interest to know whether or not
this algebra can also be made to be covariant under the Lie-Poisson
action of the translation group, and
hence under the action of the full Poincar\'e group.
With regard to Lorentz transformations alone, it is also
of interest to know whether or not the angular momenta
$\gamma$ can somehow play the role
of generators of the transformation, as in the canonical theory, and also
whether or not the momenta $\tilde p$ can somehow play the role
of generators of translations.   One thing which is clear is
that infinitesimal Lorentz
transformations are not obtained (as in the canonical theory) by
simply taking Poisson brackets with $\gamma$.
Similarly, translations are not obtained
by simply taking Poisson brackets with $\tilde p$.
[An analogous problem was solved in \cite{MSS} upon studying the
system of a isotropic rigid rotator.  That system was invariant under
the Lie-Poisson action of the chiral symmetry group.
There we were able to find
 the generators of the chiral symmetry, and
they took values in a group which was dual to the symmetry
group.  Similar novel features are anticipated here.]

In Secs. 3 and 4, realizations for the deformed Poincar\'e
algebra were found  which were associated with a
 single relativistic particle.  It is then
 of interest to know how one constructs representations for
  two or more particles.  This is not straightforward because as we
  found in Sec. 4, the Poisson structure for the observables $\gamma$
is not compatible with
the product for the Lorentz group, i.e. the space spanned by $ \{\gamma\}$
did not define a Poisson-Lie group.  Also, it can be checked that
the Poisson structure for the observables $\tilde p$
is not compatible with addition of the momenta.
Just as we found in the case of a single particle that
the spin does not commute with the orbital
angular momentum, we can conclude
the angular momenta for different particles cannot commute, and
we also suspect that the momenta
of different particles does not commute.

\bigskip
\bigskip

{\parindent 0cm{\bf Acknowledgements:}}
We are grateful to G. Marmo for helpful discussions.
We were supported in part
by the Department of Energy, USA, under contract number
DE-FG05-84ER40141.

\bigskip
\bigskip
{\parindent 0cm{\bf Appendix}}
\bigskip

Here we show that the Poisson brackets (\ref{pbpf}) for $\tilde p$
with $f$ can be deduced using the realization for $f$ in terms of
$x$ and $\tilde p$  given in eq. (\ref{feiJ}),
along with the Poisson brackets (\ref{pbptpt}) and (\ref{pbxpt}).

We start by
computing the brackets for $\tilde p$ with sin $\lambda J=\lambda
x\tilde p$.  From eqs. (\ref{pbptpt}) and (\ref{pbxpt}) we get
\begin{eqnarray}
   \{\ptone,\sJtwo\}&=&   r^\dagger \ptone\sJtwo +\ptone\sJtwo r^\dagger
  -\sJtwo r^\dagger \ptone -\ptone r^\dagger \sJtwo \cr
 & & +\;\lambda (\fone^\dagger+\ftwo)   \ptone  \Pi\;    \;.\label{pbpsJ}
\end{eqnarray}
To determine the brackets of
 $\tilde p$ with $f$, we also need  $\pois{\ptone}{\cJtwo}$.
We can deduce it by knowing
the brackets for $\tilde p$ with $\cos^2 \lambda J=1-\sin^2 \lambda J$,
  which are easily obtained from eq. (\ref{pbpsJ}),
\begin{eqnarray}
   \{\ptone,\c2Jtwo\}&=&   r^\dagger \ptone\c2Jtwo +\ptone\c2Jtwo
   r^\dagger -\c2Jtwo r^\dagger \ptone -\ptone r^\dagger \c2Jtwo \cr
& &-\;\lambda \biggl(\sJtwo(\fone^\dagger+\ftwo)   +
(\fone^\dagger+\ftwo)   \sJdone \biggr)\ptone\Pi \;  \;,\label{pbpc2J}
\end{eqnarray}
where we have used
\Eq{\tilde p\; \sin (\lambda J)=\sin (\lambda J^\dagger)\; \tilde p \;.
\label{pssp} } Then a solution is
\begin{eqnarray}
   \{\ptone,\cJtwo\}&=&   r^\dagger \ptone\cJtwo +\ptone\cJtwo
   r^\dagger -\cJtwo r^\dagger \ptone -\ptone r^\dagger \cJtwo \cr
& &-\;i\lambda (\fone^\dagger-\ftwo)\ptone\Pi \;  \;.\label{pbpcJ}
\end{eqnarray}
To check that eq. (\ref{pbpc2J}) follows from eq.
(\ref{pbpcJ}), we can apply the identities
\Eq{\tilde p\; \cos (\lambda J)=\cos (\lambda J^\dagger)\; \tilde p \;,
\label{pccp}  }  and
\Eq{  \sJtwo(\fone^\dagger+\ftwo) +
(\fone^\dagger+\ftwo)   \sJdone =
i \cJtwo(\fone^\dagger-\ftwo) + i(\fone^\dagger-\ftwo) \cJdone \;.}
The former identity follows after writing
 $J$ as a polynomial function of $x\tilde p$, while the latter
  follows after writing $f=\cos(\lambda J)+i\sin (\lambda J)$.
Finally,  from eqs. (\ref{pbpsJ}) and (\ref{pbpcJ}), and
 using a third identity, i.e. eq. (\ref{rrdid}),
we then obtain the desired result eq. (\ref{pbpf}).

\end{document}